\documentclass[11pt]{article}
\usepackage[ansinew]{inputenc}
\usepackage{array}
\usepackage{color}
\usepackage{amsmath,amsthm,amscd}
\usepackage{amsxtra}
\usepackage{amstext}
\usepackage{amssymb}
\usepackage{latexsym}
\usepackage{amsfonts,cite}
\usepackage{graphics,epstopdf}
\usepackage{epsfig}
\usepackage{amsfonts}
\usepackage{graphicx}
\usepackage{titlesec}
\titleformat{\section}
{\normalfont\Large\bfseries}
{\thesection.}{1em}{}
\usepackage{tcolorbox}
\usepackage{array,tabularx,colortbl}
\usepackage[framemethod=TikZ]{mdframed}
\usepackage{hyperref}

\numberwithin{equation}{section}
\renewenvironment{proof}{{\bfseries Proof:}}{\qed} 
\providecommand{\U}[1]{\protect\rule{.1in}{.1in}}

\begin{document}
	
	\sloppy

	\thispagestyle{empty}
	\parindent=0mm
	\begin{center}
		{\Large{\bf Bivariate degenerate Hermite polynomials in the framework of Lie algebra $\mathcal{K}_5$ }}\\
		\vspace{0.3cm}
		{\bf Subuhi Khan$^1$ and Mahammad Lal Mia$^2$}\\
		\vspace{0.2cm}
		Department of Mathematics\\
		Aligarh Muslim University\\
		Aligarh, India\\
		\footnote{Emails:}
		\footnote{1. subuhi2006@gmail.com (Subuhi Khan)} 
		\footnote{2. lalmia140511@gmail.com(Mahammad Lal Mia)}
		\footnote{\\
			This work has been done under  Junior Research Fellowship (Ref. No. 221610098808 dated:28/10/2022) awarded to the second author by the University  Grants Commission, Government of India, New Delhi.}
	\end{center}
	\parindent=0mm
	\vspace{0.1cm}
	\noindent
	{\bf Abstract:}~
	In this article, the matrix elements of the representation $\uparrow'_{\omega,\mu}$ of the 5-dimensional Lie algebra $\mathcal{K}_5$ are obtained for the first time. The bivariate degenerate Hermite polynomials $\mathcal{H}_m(z_1,z_2|\tau)$ are considered within the context of this representation. Further, employing the Lie algebraic techniques, certain specific results concerning these polynomials are established. Some examples providing the implicit formulas for the polynomials related to the polynomials $\mathcal{H}_m(z_1,z_2|\tau)$ are considered. Integral equations for these polynomials are also explored.\\
	
	
	
	\parindent=8mm
	
	\section{Introduction} 	There are numerous approaches to comprehend the usage of classical special functions in mathematical physics. These approaches include studying orthogonal polynomials, solving second-order ordinary differential equations, using integral transforms of other functions, deriving expressions by using contour integrals in the complex plane, specializing hypergeometric or related series, and exploring bonds with Lie groups. Weisner \cite{Weisner} employed a Lie theoretic approach to derive generating functions that underscore the importance of group theory in relation to Hermite, Bessel, hypergeometric and certain other special functions. Miller \cite{Miller, MillerW} systematically presented a detailed study of Lie-theoretic approach to special functions, establishing a solid foundation for it. 
	The connections between generalized special functions and Lie group theory has been recognized by some researchers, see for example,  \cite{ sub, subu, subuhi}.\\
	
	\noindent Special functions serve as basis vectors or matrix elements for the unitary irreducible representations of Lie groups of low-dimensions \cite{James, Vilen}. Investigating these Lie groups along with corresponding Lie algebras provides a cohesive approach to a substantial portion of special functions theory, particularly the aspects applicable to mathematical physics. \\
	
	\noindent The main tools necessary for deriving the desired results include multiplier representations of local Lie groups and representations of Lie algebras via generalized Lie derivatives. For composing irreducible representations, there is a regular basis, including Bessel, hypergeometric, and confluent hypergeometric functions. This relationship between special functions and Lie algebras offers valuable insights into special functions theory.\\
	
	\noindent $K_5$ is the complex Lie group of dimension five and general component of this complex Lie group is $g(q, \alpha, \beta, \gamma, \delta)$, \quad $q, \alpha, \beta, \gamma, \delta \in \mathbb{C}$. The group has following $5\times5$ matrix realization: 
	\begin{equation} \label{hpeq1}
		g(q, \alpha, \beta, \gamma, \delta) = \begin{pmatrix}
			1 & \gamma e^\delta& \beta e^{-\delta}& 2\alpha-\beta\gamma &\delta \\
			0& e^\delta& 2qe^{-\delta}& \beta-2q\gamma& 0 \\
			0& 0& e^{-\delta}& -\gamma& 0 \\
			0&0&0&1&0 \\
			0&0&0&0&1
		\end{pmatrix},\quad \alpha, \beta, \gamma, \delta, q  \in \mathbb{C}.
	\end{equation}\\
\noindent	The group operation performed through matrix operation and multiplicative law of this group is 
	\begin{equation*}
		g(q, \alpha, \beta, \gamma, \delta) g(q', \alpha', \beta', \gamma', \delta') = g(q+e^{2\delta}, \alpha+\alpha'+e^\delta \gamma \beta', \beta+e^\delta \beta'+2e^{2\delta} \gamma q', \gamma+e^{-\delta} \gamma', \delta+\delta').
	\end{equation*}
	It is to be noted that the identity element of $K_5$ is $g(0, 0, 0, 0, 0)$,  whereas the inverse of $g(q, \alpha, \beta, \gamma, \delta)$ is $g(-qe^{-2\delta}, -a+\beta\gamma-2\gamma^2 q, -\beta e^-\delta +2e^{-\delta} \gamma q', -\gamma e^{\delta} , -\delta)$.
	All the elements of this group set up a subgroup of $K_5$ for $q=0$, which is isomorphic to the harmonic oscillator Lie group $G(0,1)$ \cite[p. 9]{Miller}.\\
	
	\noindent The corresponding Lie algebra $\mathcal{K}_5$ possesses the elements of following type: \\
	\begin{equation*}
		\begin{pmatrix}
			0& \eta_4&\eta_3  &\eta_2&\eta_5 \\
			0& \eta_5& 2\eta_1& \eta_3& 0 \\
			0& 0& -\eta_5& -\eta_4& 0 \\
			0&0&0&1&0 \\
			0&0&0&0&1
		\end{pmatrix},
	\end{equation*}
	
\noindent	where $\eta_1=\frac{dq}{dt}$, $\eta_2=\frac{d\alpha}{dt}$, $\eta_3=\frac{d\beta}{dt}$, $\eta_4=\frac{dc}{dt}$ and $\eta_5=\frac{d \delta}{dt}$ in the neighborhood of identity.\\ 
	
	\noindent For $q=0$ the Lie algebra 
	$\mathcal{K}_5$ turns into the Lie algebra $L(G(0,1))$ of $G(0,1)$, denoted by $\mathcal{G}(0,1)$. The basis elements of $\mathcal{K}_5$ are\\
	
		\noindent  
	$\mathsf{j}^{+} = \begin{pmatrix}
		0& 0&1 &0&0 \\
		0& 0& 0& 1& 0 \\
		0& 0&0 & 0& 0 \\
		0&0&0&0&0 \\
		0&0&0&0&0
			\end{pmatrix}$, \quad $\mathsf{j}^{-} = \begin{pmatrix}
		0& 1&0 &0&0 \\
		0& 0& 0& 0& 0 \\
		0& 0&0 & -1& 0 \\
		0&0&0&0&0 \\
		0&0&0&0&0
		\end{pmatrix}$, \quad $\mathsf{j}^{3} = \begin{pmatrix}
		0& 0&0&0&1 \\
		0& 1& 0& 1& 0 \\
		0& 0&-1& 0& 0 \\
		0&0&0&0&0 \\
		0&0&0&0&0
			\end{pmatrix}$, \\
	
\hspace{-0.8cm}	$\mathcal{E} = \begin{pmatrix}
		0& 0&0&1&0 \\
		0& 0& 0& 0& 0 \\
		0& 0&0 & 0& 0 \\
		0&0&0&0&0 \\
		0&0&0&0&0
			\end{pmatrix}$, \quad $\mathcal{Q} = \begin{pmatrix}
		0& 0&0 &0&0 \\
		0& 0& 0& 1& 0 \\
		0& 0&0 & 0& 0 \\
		0&0&0&0&0 \\
		0&0&0&0&0
		\end{pmatrix}$,\\ 
	
	\noindent which satisfy the following commutation relations:\\
	\begin{equation*}
		[\mathsf{j}^{3}, \mathcal{Q}] = 2 \mathcal{Q}, \quad
		[\mathsf{j}^3, \mathsf{j}^{\pm}] = \pm \mathsf{j}^{\pm}, 
	\end{equation*}
	\begin{equation} \label{hpeq2}
		[\mathsf{j}^{-}, \mathcal{Q}] = 2 \mathsf{j}^{+}, \quad [\mathsf{j}^{-}, \mathsf{j}^{+}] = \mathcal{E}, \quad [\mathsf{j}^{+}, \mathcal{Q}] = \mathcal{O}, 
	\end{equation}
	\begin{equation*}
		[\mathsf{j}^{3},\mathcal{E}] = [\mathcal{Q}, \mathcal{E}] = [\mathsf{j}^{\pm}, \mathcal{E}] =     \mathcal{O}, 
	\end{equation*}
	where $\mathcal{O}$ is the $5\times5$ zero matrix.\\
	
	\noindent The Hermite polynomials $\mathcal{H}_m(z)$ \cite{Andrews} and Hermite numbers $\mathcal{H}_m$ are defined by the following relations:
	\begin{equation*} 
		e^{2zt-t^2}=\sum_{m=0}^{\infty}\mathcal{H}_m(z)\frac{t^m}{m!}
	\end{equation*}
	and
	\begin{equation*}
		e^{-t^2}=\sum_{m=0}^{\infty}\mathcal{H}_m\frac{t^m}{m!},
	\end{equation*}
	respectively.\\
	The series expansion of Hermite polynomials is given as
	\begin{equation} \label{hpeq83}
		\mathcal{H}_m(z)=m! \sum_{s=0}^{\left[ \dfrac{m}{2} \right]} \dfrac{(-1)^s (2z)^{m-2s}}{s! (m-2s)!}.
	\end{equation}
		\noindent Hermite numbers and polynomials hold significant importance in various fields, including number theory, combinatorics, special functions, and differential equations. These polynomials are a classical orthogonal polynomial sequence, manifesting within the Edgeworth series in probability theory. In physics, these polynomials correspond to the eigen states of the quantum harmonic oscillator.\\
	
	\noindent
	Special functions of two variables, like spherical harmonics and hypergeometric functions, are crucial for solving partial differential equations and representing physical systems in various coordinate frameworks. Their unique properties, including orthogonality and roles as eigen functions, make them indispensable in fields such as quantum mechanics, fluid dynamics, and computational analysis, where they simplify complex mathematical and engineering problems. The 2-variable Hermite polynomials $\mathcal{H}_m(z_1,z_2)$ possess the generating  relation \cite{AppHer}
	\begin{equation} \label{hpeq3}
		e^{z_1t+z_2t^2} = \sum_{m=0}^{\infty} \mathcal{H}_m(z_1,z_2) \frac{t^m}{m!},
	\end{equation} 
	whereas the series expansion of $\mathcal{H}_m(z_1,z_2)$ is given as
	\begin{equation} \label{hpeq4}
		\mathcal{H}_m(z_1,z_2) = m! \sum_{s=0}^{\left[ \dfrac{m}{2} \right]} \frac{z_1^{m-2s} z_2^s}{(m-2s)! s!}.
	\end{equation}
	\noindent
	The significance and applications of these polynomials in various fields have led researchers to study various special polynomials and numbers, with particular attention to degenerate forms of Stirling, Bernoulli, Euler, Genocchi, and tangent polynomials \cite{Carlitz, CenHow}.\\
	
	\noindent The 2-variable degenerate Hermite polynomials (2VDHP) $\mathcal{H}_{m}(z_1,z_2|\tau)$ are introduced in \cite{Hwang} with the following generating function:
	\begin{equation} \label{hpeq5}
		\left(1+\tau \right)^{\dfrac{t(z_1+z_2t)}{\tau}}=\sum_{m=0}^{\infty}\mathcal{H}_{m}(z_1,z_2|\tau) \frac{t^m}{m!}.
	\end{equation}
	\noindent If $\tau \rightarrow 0$, then $(1+\tau)^{\frac{t}{\tau}} \rightarrow e^t$. Therefore, for $\tau \rightarrow 0$ , it is clearly seen that generating equation \eqref{hpeq5} reduces to generating function \eqref{hpeq3}. The series expansion of 2VDHP $\mathcal{H}_{m}(z_1,z_2|\tau)$ is given by
	\begin{equation} \label{hpeq6}
		\mathcal{H}_{m}(z_1,z_2|\tau) = \sum_{s=0}^{\left[ \dfrac{m}{2} \right]} \frac{m!}{s! (m-2s)!}\left(\frac{\log(1+\tau)}{\tau}\right)^{m-s}z_1^{m-2s} z_2^s  .
	\end{equation}\\
	\noindent The recursive relations and differential equation satisfied by 2VDGHP $\mathcal{H}_{m}(z_1,z_2|\tau)$, are derived by proving the following Lemma:\\
	
		\noindent	{\bf Lemma 1.1}
	For any positive integer $m$, the following recursive relations and differential equation for the 2-variable degenerate Hermite polynomials $\mathcal{H}_{m}(z_1,z_2|\tau)$ hold:
	
	\begin{equation} \label{hpeq9}
		\dfrac{\partial}{\partial z_1}	\mathcal{H}_{m}(z_1,z_2|\tau)=m\frac{\log(1+\tau)}{\tau}\mathcal{H}_{m-1}(z_1,z_2|\tau), 
	\end{equation}	
	\begin{equation} \label{hpeq10}
		\dfrac{\partial^s}{\partial z_1^s}\mathcal{H}_m(z_1,z_2|\tau)=\frac{m!}{(m-s)!} \left({\frac{\log(1+\tau)}{\tau}}\right)^s \mathcal{H}_{m-s}(z_1,z_2|\tau), \quad s\leq m, 
	\end{equation}
	\begin{equation} \label{hpeq11}
		\dfrac{\partial}{\partial z_2}	\mathcal{H}_{m}(z_1,z_2|\tau)=m(m-1)\frac{\log(1+\tau)}{\tau}\mathcal{H}_{m-2}(z_1,z_2|\tau), 
	\end{equation}
	\begin{equation} \label{hpeq12}
		\dfrac{\partial^s}{\partial z_2^s}\mathcal{H}_m(z_1,z_2|\tau)=\frac{m!}{(m-2s)!} \left({\frac{\log(1+\tau)}{\tau}}\right)^s \mathcal{H}_{m-2s}(z_1,z_2|\tau), \quad s\leq \left[\frac{m}{2}\right], 
	\end{equation}	
	\begin{equation} \label{hpeq13}
		\mathcal{H}_{m+1}(z_1,z_2|\tau)= \frac{\log(1+\tau)}{\tau}z_1\mathcal{H}_m(z_1,z_2|\tau)+2z_2m\frac{\log(1+\tau)}{\tau}\mathcal{H}_{m-1}(z_1,z_2|\tau), 
	\end{equation}	
	and
	\begin{equation} \label{hpeq14}
		\left(2z_2\dfrac{\partial^2}{\partial z_1^2} + z_1\frac{\log(1+\tau)}{\tau}\dfrac{\partial}{\partial z_1} -m\frac{\log(1+\tau)}{\tau}\right)\mathcal{H}_{m}(z_1,z_2|\tau)=0. 
	\end{equation} 
	\noindent	\begin{proof}
		Differentiating generating equation \eqref{hpeq5} w.r.t. $z_1, z_2,$ and $t$, we obtain differential and pure recurrence relations \eqref{hpeq9}, \eqref{hpeq11} and \eqref{hpeq13}, respectively. Again, differentiating equation \eqref{hpeq5}, $s$-$th$ times with respect to $z_1$ and $z_2$, relations \eqref{hpeq10} and \eqref{hpeq12} are proved. Further, replacing $m$ by $m-1$ in equation \eqref{hpeq13} and using relations \eqref{hpeq10} and \eqref{hpeq12} (for $s=2$), we find differential equation \eqref{hpeq14}.
	\end{proof}\\
	
	\noindent In this article, the 2-variable degenerate Hermite polynomials 2VDHP $\mathcal{H}_{m}(z_1,z_2|\tau)$ are encased into the surrounding of an irreducible representation $\uparrow'_{\omega,\mu}$ of $\mathcal{K}_5$. The article is organized in four sections. In section 2, the general matrix elements of the representation $\uparrow'_{\omega,\mu}$ for the Lie algebra $\mathcal{K}_5$ are computed. Further, the implicit formulas involving the 2VDHP $\mathcal{H}_{m}(z_1,z_2|\tau)$ are established using the matrix elements corresponding to this representation and also for certain specific parameter values. Section 3, explores some examples to derive implicit formulas for the polynomials related to the 2VDHP $\mathcal{H}_{m}(z_1,z_2,\tau)$. In section 4, the Volterra integral equation for these polynomials is derived.
	\section{Representation $\uparrow'_{\omega,\mu}$ of $\mathcal{K}_5$ and $\mathcal{H}_{m}(z_1,z_2|\tau)$}
		
	\noindent Let $\rho$ be an abstract irreducible representation of $\mathcal{K}_5$ on a vector space $V$, which satisfies the following property:\\
	
	\noindent	$\Large \mathcal{\boldmath P}:$ $\rho/_{\mathcal{G}(0,1)}$ ($\rho$ restricted to the sub-algebra $\mathcal{G}(0,1)$ ) is isomorphic to the analogues irreducible representations of $\mathcal{G}(0,1)$.\\
	and let 
	\begin{equation*}
		\rho(\mathsf{j}^{\pm})=J^{\pm}, \quad 	\rho(\mathsf{j}^{3})=J^{3},\quad \rho(\mathcal{E})=E, \quad   \rho(\mathcal{Q})=Q
	\end{equation*}
	 be operators on $V$.
	Operators $ J^{\pm}, J^3, E, Q$ hold relations identical to \eqref{hpeq2}.
	Here, the representation $\rho$ of $\mathcal{K}_5$ satisfying property $\Large \mathcal{\boldmath P}$ is considered. This representation is isomorphic to the representation $\uparrow'_{\omega,\mu}$, \quad $\omega$, $\mu$ $\in \mathbb{C}$, with $\mu \neq 0$. The set of eigen values of $J^3$ is the set $\mathcal{S} = \{-\omega +r : r ~is~ a~ nonnegative~ integer\}.$ The representation space $V$ has a basis $\{f_m\}$, $m \in \mathcal{S}$ with following relations: 
	\begin{equation*}
		J^+f_m = \mu f_{m+1}, \quad J^-f_ = (\omega+r)f_{m-1},\quad J^3f_m = mf_m 
	\end{equation*}
	\begin{equation} \label{hpeq7}
		Ef_m = \mu f_m,  \quad Qf_m = \mu f_{m+2}.
	\end{equation}
	
	\noindent A realization of the representation $\uparrow'_{\omega,\mu}$ given by some $\tau \in \mathcal{O}(\tilde{\alpha}_1)$, set or cohomology class of all realizations of the representation that are equivalent (up to cohomology), acting on a vector space $\Psi$ of analytic functions of $z$ is considered in  \cite{Miller}. A realizations of irreducible representation $\uparrow'_{\omega,\mu}$ of $\mathcal{K}_5$ is constructed in such a way that the operators $ J^{\pm}, J^3, E, Q$ are given by
	\begin{equation}
		J^3=z\dfrac{d}{dz}-\omega,\quad J^+=\mu z, \quad J^-=\dfrac{d}{dz},\quad E=\mu, \quad Q=\mu z^2.
	\end{equation}
	 The space of all finite linear combinations of the functions $h_r(z)=z^{r}, r=0,1,2, \ldots$ is designated by $\Psi_2$ and the basis vectors $f_m$ of $\Psi_2$ are defined by $f_m(z)=h_r(z)=z^{r}$ where $m=-\omega +r$, $r \geq 0$. Evidently
	\begin{align}
		&J^3f_m=\left(z\dfrac{d}{dz}-\omega\right)z^r=mf_m,\quad J^+f_m=(\mu z)z^r=\mu f_{m+1}, \quad J^-f_m=\dfrac{d}{dz}z^r=(m+\omega)f_{m-1}, \nonumber \\
		&\hspace{4cm} Ef_m=\mu f_m, \quad Qf_m=\left(\mu z^2\right)z^r=\mu f_{m+2}.
	\end{align}
	\noindent These relations describe a realizations of $\uparrow'_{\omega,\mu}$ of $\mathcal{K}_5$ on $\Psi_2$. This realization is extended local multiplier representation $B$ of $K_5$ on the space $\mathcal{O}_2$ consisting of entire functions of $z$. The operators $B(g)$, $g=g(q,\alpha,\beta,\gamma,\delta ) \in K_5$ defining this representation is given by
	\begin{equation} \label{hpeq150}
		[B(g)f](z)=\exp[\mu (qz^2+\beta z+\alpha)-\omega \delta]f(e^\delta z+e^\delta \gamma), \quad  f \in \mathcal{O}_2.
	\end{equation}
	\noindent The matrix elements $A_{lr}(g)$ of the operators $B(g)$ with respect to the basis $\{f_m=h_r\}$ are given by
	\begin{equation*}
		[B(g)h_r](z)=\sum_{l=0}^{\infty}A_{lr}(g)h_l(z), \quad r=0,1,2 \ldots, 
	\end{equation*}
	which in view of representation \eqref{hpeq150}, takes the form
	\begin{equation}\label{hpeq151}
		\exp[\mu (qz^2+\beta z+\alpha)+(r-\omega) \delta](z+\gamma)^r=\sum_{l=0}^{\infty}A_{lr}(g)h_l(z).
	\end{equation}
	Miller \cite{Miller} computed matrix elements corresponding to some special choices of the group parameters. However, the expression for general matrix elements is not computed in \cite{Miller}. In order to determine the expressions for the general matrix elements of the representations $\uparrow'_{\omega,\mu}$ for Lie algebra $\mathcal{K}_5$ with respect to the basis $f_m(z)=h_r(z)=z^{r}$, where $m=-\omega +r$, $r \geq 0$, the following result is proved:\\
	
	\noindent	{\bf Theorem 2.1}
		 The matrix elements $A_{lr}(g)$ of the representations $\uparrow'_{\omega,\mu}$ for the Lie algebra $\mathcal{K}_5$ with respect to basis $f_m(z)=h_r(z)=z^{r}$, where $m=-\omega +r$, $r \geq 0$ are given by 
	\begin{equation} \label{hpeq17}
		A_{lr}(g)=\exp(\mu \alpha +(r-\omega)\delta) \frac{\gamma^{r-l} r!}{(r-l)! l!}F_{1\colon 1; 1}^{1\colon 1; 1} \left( \begin{matrix}
			[(-l):2,1]: [(1):1];[(1):1];&\\
			& \mu q \gamma ^2,-\mu \beta \gamma \\
			[(r-l+1):2,1]: [(1):1];[(1):1];
		\end{matrix} \right),
	\end{equation}  	
	\hspace{13cm} $r\geq l\geq 0$.\\
	or 
	\begin{equation} \label{hpeq85}
		A'_{lr}(g)=\exp[\mu \alpha+(r-\omega)\delta]  \sum_{j=0}^{\infty}\frac{\gamma^{j} r! \left(-\mu q \right)^{\frac{l-r+j}{2}} H_{l-r+j}\left(\frac{\mu \beta }{2(-\mu q)^{\frac{1}{2}}} \right) }{(r-j)! j!(l-r+j)!}, \quad  l\geq r\geq 0. 
	\end{equation}

\noindent	\begin{proof}
		Expanding the exponentials containing powers of $z$ of equation \eqref{hpeq151} and simplifying, it follows that
	\begin{equation} \label{hpeq19}
		\exp[\mu \alpha +(r-\omega)\delta ]\sum_{j=0}^{r}  \sum_{s,k=0}^{\infty}\frac{\gamma^{r-j}r!(\mu q)^s (\mu \beta)^k}{(r-j)! j! s! k!} z^{2s+k+j} =\sum_{l=0}^{\infty}A_{lr}(g)z^l. 
	\end{equation}
	Replacing $2s+k+j$ by $l$ and rearranging the series, we have
	\begin{equation} \label{hpeq20}
		\exp[\mu \alpha+(r-\omega)\delta]\sum_{l=0}^{\infty}  \sum_{s,k=0}^{\infty}\frac{\gamma^{r-l}r!\left(\mu q \gamma ^2\right)^s (\mu \beta \gamma)^k}{(r-l+2s+k)! (l-2s-k)! s! k!} z^l =\sum_{l=0}^{\infty}A_{lr}(g)z^l, 
	\end{equation}
	which on equating the like powers of $z$, gives
	\begin{equation} \label{hpeq21}
		A_{lr}(g)=\exp[\mu \alpha+(r-\omega)\delta]  \sum_{s,k=0}^{\infty}\frac{\gamma^{r-l} r! \left(\mu q \gamma^2\right)^s (\mu \beta \gamma)^k}{(r-l+2s+k)!(l-2s-k)!s!k!}. 
	\end{equation}
	Next, Let us recall that the generalized Lauricella function of two variables is defined by the following series \cite[p. 64(18)]{SRIV}:
	\begin{align} \label{hpeq22}
		F_{C\colon D; D'}^{A\colon B; B'} \left( \begin{matrix}
			z_1\\
			z_2
		\end{matrix} \right) \equiv	&F_{C\colon D; D'}^{A\colon B; B'}\left( \begin{matrix}
			[(a):v,\phi]: [(b):\psi];[(b'):\psi'];&\\
			& z_1,z_2\\
			[(c):\delta,\epsilon]: [(d):\eta];[(d'):\eta'];
		\end{matrix} \right) \nonumber \\
		&=\sum_{s,k=0}^{\infty} \frac{ \displaystyle \prod_{j=1}^{A}(a_j)_{sv_j+k\phi_j} \displaystyle \prod_{j=1}^{B}(b_j)_{s\psi_j} \displaystyle \prod_{j=1}^{B'}(b'_j)_{k\psi'_j} z_1^s z_2^k}{\displaystyle \prod_{j=1}^{C}(c_j)_{s \delta_j+k\epsilon_j} \displaystyle \prod_{j=1}^{D}(d_j)_{s\eta_j} \displaystyle \prod_{j=1}^{D'}(d'_j)_{k\eta'_j}s!k!},
	\end{align}  	
	which converges absolutely, when
	\begin{equation} \label{hpeq23}
		\Delta =1+ \sum_{j=1}^{C} \delta_j + \sum_{j=1}^{D} \eta_j-\sum_{j=1}^{A} v_j - \sum_{j=1}^{B} \psi_j>0;\quad \Delta' =1+ \sum_{j=1}^{C} \epsilon_j + \sum_{j=1}^{D'} \eta'_j-\sum_{j=1}^{A} \phi_j - \sum_{j=1}^{B'} \psi'_j>0. \nonumber
	\end{equation} 
	For $A=B=B'=C=D=D'=1$, $a_1=-l, c_1=r-l+1, v_1=\delta_1=2$, $\phi_1=\psi_1=\psi'_1=\epsilon_1=\eta_1=\eta'_1=1$, and $b_1=d_1=b'_1=d'_1=1$, series \eqref{hpeq22} remains convergent and takes the following form:
	\begin{equation} \label{hpeq24}
		F_{1\colon 1; 1}^{1\colon 1; 1} \left( \begin{matrix}
			[(-l):2,1]: [(1):1];[(1):1];&\\
			& z_1,z_2\\
			[(r-l+1):2,1]: [(1):1];[(1):1];
		\end{matrix} \right) =\sum_{s,k=0}^{\infty} \frac{(-l)_{2s+k} (1)_s (1)_k}{(r-l+1)_{2s+k}(1)_s (1)_k}\frac{z_1^s}{s!}\frac{z_2^k}{k!},
	\end{equation}
	which on simplifications, becomes
	\begin{align} \label{hpeq65}
		\notag
		&\hspace{-0.55cm}	F_{1\colon 1; 1}^{1\colon 1; 1} \left( \begin{matrix}
			[(-l):2,1]: [(1):1];[(1):1];&\\
			& z_1,z_2\\
			[(r-l+1):2,1]: [(1):1];[(1):1];
		\end{matrix} \right)=\sum_{s,k=0}^{\infty} \frac{(l)! (r-l)!}{(r-l+2s+k)!(l-2s-k)!}\\
		&\hspace{12cm}\times\frac{z_1^s}{s!}\frac{(-z_2)^k}{k!}.
	\end{align}
	
	\noindent Finally making use of series definition \eqref{hpeq65} in equation \eqref{hpeq21},  expressions \eqref{hpeq17} for matrix elements $A_{lr}(g)$ are obtained.\\
	
	\noindent Again, considering an alternate series expansion of equation \eqref{hpeq151} in the l.h.s. and denoting the corresponding matrix elements in the r.h.s. by $A'_{lr}(g)$, it follows that
	\begin{equation} \label{hpeq80}
		\exp[\mu \alpha +(r-\omega)\delta ]\sum_{j=0}^{\infty}  \sum_{s,k=0}^{\infty}\frac{\gamma^{j}r!(\mu q)^s (\mu \beta)^k}{(r-j)! j! s! k!} z^{2s+k+r-j} =\sum_{l=0}^{\infty}A'_{lr}(g)z^l. 
	\end{equation}
	Replacing $2s+k+r-j$ by $l$ and rearranging the series, we have
	\begin{equation} \label{hpeq81}
		\exp[\mu \alpha+(r-\omega)\delta]\sum_{l=0}^{\infty}  \sum_{j=0}^{\infty}\sum_{s=0}^{\scalebox{0.6}{$\left[ \dfrac{l - r + j}{2} \right]$}}\frac{\gamma^{j}r!\left(\mu q \right)^s (\mu \beta)^{l-r+j-2s}}{(l-r+j-2s)! (r-j)! j! s!} z^l =\sum_{l=0}^{\infty}A'_{lr}(g)z^l, 
	\end{equation}
	which on equating the like powers of $z$, yields
	\begin{equation} \label{hpeq82}
		A'_{lr}(g)=\exp[\mu \alpha+(r-\omega)\delta]  \sum_{j=0}^{\infty}\sum_{s=0}^{\scalebox{0.6}{$\left[ \dfrac{l - r + j}{2} \right]$}}\frac{\gamma^{j} r! \left(\mu q \right)^s (\mu \beta )^{l-r+j-2s}}{(l-r+j-2s)!(r-j)!j!s!}. 
	\end{equation}
	or, equivalently
	\begin{equation} \label{hpeq84}
		A'_{lr}(g)=\exp[\mu \alpha+(r-\omega)\delta]  \sum_{j=0}^{\infty}\sum_{s=0}^{\scalebox{0.6}{$\left[ \dfrac{l - r + j}{2} \right]$}}
		\frac{ (-1)^s \gamma^{j} r! \left(-\mu q \right)^{\frac{l-r+j}{2}} \left(\frac{\mu \beta }{(-\mu q)^{\frac{1}{2}}} \right)^{l-r+j-2s} }{(l-r+j-2s)!(r-j)!j!s!}. 
	\end{equation}
	Use of definition \eqref{hpeq83} in expansion \eqref{hpeq84}, yields expression \eqref{hpeq85} for matrix elements $A'_{lr}(g)$.\\
\end{proof} \\
	\noindent	{\bf Remark 2.1}
		 Taking $q=0$ in expressions \eqref{hpeq17} and \eqref{hpeq85} and denoting the corresponding matrix elements by $B_{lr}(g)$ and $B'_{lr}(g)$ respectively, we find
		 \begin{equation*}
		 	B_{lr}(g)=\exp(\mu \alpha +(r-\omega)\delta) \frac{\gamma^{r-l} \quad \Gamma(r+1)}{(r-l)! \quad \Gamma(l+1)} {}_1F_1
		 	(-l;r-l+1;-\mu \beta \gamma), \quad \text{if} \quad r \geq l\geq 0,
		 \end{equation*}
		 	\begin{equation*}
		 		B'_{lr}(g)=\exp(\mu \alpha +(r-\omega)\delta) \frac{(\mu \beta)^{l-r}}{(l-r)!} {}_1F_1
		 		(-r;l-r+1;-\mu \beta \gamma), \quad \text{if} \quad l \geq r\geq 0.
		 	\end{equation*}
		 	Here ${}_1F_1$ denotes the confluent hypergeometric function \cite{Andrews}. The most convenient form to express the matrix elements in terms of the generalized Laguerre polynomials ${L_{n}}^{\alpha}(x)$ \cite{Andrews} is given as 
\begin{equation} \label{hpeq66}
B_{lr}(g)	=\exp(\mu \alpha+(r-\omega)\delta)\gamma^{r-l}{L_{l}}^{r-l}(-\mu \beta\gamma), \quad k,l\geq0.
\end{equation}
	\noindent	{\bf Remark 2.2}
		Taking $\alpha=\beta=\delta=0$ in expressions \eqref{hpeq17} and \eqref{hpeq85} and  denoting the corresponding matrix elements by $C_{lr}(g)$, we find 
	\begin{equation} \label{hpeq68}
		C_{lr}(g)=\gamma^{r-l} r!\sum_{j}  \frac{(\mu q\gamma^2)^j}{j! (2j+r-l)!(l-2j)},
	\end{equation}
	where $j$ ranges over all integral values such that the summand is defined.\\
	
\noindent	{\bf Remark 2.3}
		Taking $\alpha=\gamma=\delta=0$ in expressions \eqref{hpeq17} and \eqref{hpeq85} and denoting the resultant matrix elements by $D_{lr}(g)$, we find
	\begin{equation} \label{hpeq69}
			D_{lr}(g)=  \begin{cases}
				0 &   , if~ 0\leq l \leq r,\\
				\frac{(-\mu q)^{\frac{l-r}{2}}}{(l-r)!}\mathcal{H}_{l-r}\left(\frac{\mu \beta}{2(-\mu q)^{\frac{1}{2}}}\right) &, if ~l \geq r \geq 0.
			 \end{cases}
		 \end{equation}\\
		 	\noindent Next, we proceed to create the representation $\uparrow'_{\omega,\mu}$ of $\mathcal{K}_5$ by generalized Lie derivations on a space of two complex variables. In order to establish the results involving the 2VDHP $\mathcal{H}_{m}(z_1,z_2|\tau)$, we recognize the differential operators, whose eigen functions are expressed in the following form:
	\begin{equation} \label{hpeq8}
		f_m(z_1,z_2|\tau;t) = \mathcal{H}_{m}(z_1,z_2|\tau) t^m.
	\end{equation}
	\noindent By making use of recurrence relations \eqref{hpeq9}, \eqref{hpeq11} and \eqref{hpeq13}, the following first order linearly independent differential operators are obtained:
	\begin{align} \label{hpeq15} 
		J^+&=2z_2t \dfrac{\partial}{\partial z_1} +\frac{\log(1+\tau)}{\tau}z_1t, \nonumber\\
		J^-&=\frac{1}{t}\frac{\tau}{log(1+\tau)} \dfrac{\partial}{\partial z_1}, \nonumber\\
		J^3&=t \dfrac{\partial}{\partial t},\\
		E&=1, \nonumber\\
		Q&=2z_1z_2t^2\frac{\log(1+\tau)}{\tau} \dfrac{\partial}{\partial z_1} + 2z_2t^3 \frac{\log(1+\tau)}{\tau}\dfrac{\partial }{\partial t} + \frac{\log(1+\tau)}{\tau} \left( \frac{\log(1+\tau)}{\tau}z_1^2+2z_2\right)t^2. \nonumber
	\end{align}
	The above operators are computed in such a way that the actions of $J^+, J^-, J^3, E,$ and $Q$ on the functions $f_m(z_1,z_2|\tau;t) = \mathcal{H}_{m}(z_1,z_2|\tau) t^m$ are given by
	\begin{align} \label{hpeq16}
		J^+[\mathcal{H}_{m}(z_1,z_2|\tau) t^m]&=\mathcal{H}_{m+1}(z_1,z_2|\tau) t^{m+1},\nonumber\\
		J^-[\mathcal{H}_{m}(z_1,z_2|\tau) t^m]&= m\mathcal{H}_{m-1}(z_1,z_2|\tau) t^{m-1}, \nonumber\\
		J^3[\mathcal{H}_{m}(z_1,z_2|\tau) t^m]&= m\mathcal{H}_{m}(z_1,z_2|\tau) t^{m},\\
		E[\mathcal{H}_{m}(z_1,z_2|\tau) t^m]&=\mathcal{H}_{m}(z_1,z_2|\tau) t^{m},\nonumber\\
		Q[\mathcal{H}_{m}(z_1,z_2|\tau) t^m]&=\mathcal{H}_{m+2}(z_1,z_2|\tau) t^{m+2}.\nonumber
	\end{align}\\
		It is readily observed that $J^+, J^-, J^3, E, Q$ satisfy commutation relations identical with \eqref{hpeq2}. There is no loss of generality, in taking $\omega=0$ and $\mu =1$ for special functions theory. Thus these operators satisfy relations \eqref{hpeq7} for $\omega=0$ and $\mu =1$.\\
	
	\noindent The above observations confirm that $f_m(z_1,z_2|\tau;t) = \mathcal{H}_{m}(z_1,z_2|\tau) t^m$, set up a basis for a realization of the representation $\uparrow'_{0,1}$ of $\mathcal{K}_5$ on the space $\mathcal{F}$, which is a collection of all analytic functions in a neighborhood of the point $(z_1^0, z_2^0, t^0) = (1, 1, 1)$. The implicit formula involving the 2VDHP $\mathcal{H}_{m}(z_1,z_2|\tau)$ is established by proving the following result:\\
	
	\noindent	{\bf Theorem 2.2}
		 For the 2-variable degenerate Hermite polynomials $\mathcal{H}_{r}(z_1,z_2|\tau)$, the following implicit summation formulas holds true:
	\begin{align} \label{hpeq25} 
		&\left(1-4\xi(\tau)z_2qt^2\right)^{-\frac{r+1}{2}}\exp\left(\xi (\tau)\left(\beta^2z_2t^2+\beta z_1t\right)+\frac{\xi(\tau)}{4 z_2}z_1^2 
		\left(\left(1-4\xi(\tau)z_2qt^2\right)^{-1}-1\right)\right) \nonumber \\ &\times \mathcal{H}_r\left(\left(1-4\xi(\tau)z_2qt^2\right)^{-\frac{1}{2}} \bigg(z_1+2z_2\beta t+ 
		\dfrac{\gamma}{t}\left(1-4\xi(\tau)z_2qt^2\right)\bigg),
		z_2|\tau\right)=\sum_{l=0}^{\infty} \frac{\gamma^{r-l} r!}{l! (r-l)!}  \\
		& \times F_{1\colon 1; 1}^{1\colon 1; 1} \left( \begin{matrix}
			[(-l):2,1]: [(1):1];[(1):1];&\\
			& q \gamma ^2,-\beta \gamma \\
			[(r-l+1):2,1]: [(1):1];[(1):1];
		\end{matrix} \right)\mathcal{H}_{l}(z_1,z_2|\tau) t^{l-r}, \quad r\geq l, \nonumber
	\end{align}
	and
		\begin{align} \label{hpeq160} 
		&\left(1-4\xi(\tau)z_2qt^2\right)^{-\frac{r+1}{2}}\exp\left(\xi (\tau)\left(\beta^2z_2t^2+\beta z_1t\right)+\frac{\xi(\tau)}{4 z_2}z_1^2 
		\left(\left(1-4\xi(\tau)z_2qt^2\right)^{-1}-1\right)\right) \nonumber \\ 
		&\hspace{1cm}\times \mathcal{H}_r\left(\left(1-4\xi(\tau)z_2qt^2\right)^{-\frac{1}{2}} \bigg(z_1+2z_2\beta t+ 
		\dfrac{\gamma}{t}\left(1-4\xi(\tau)z_2qt^2\right)\bigg),
		z_2|\tau\right)  \\
		&\hspace{2cm}=\sum_{l=0}^{\infty} \sum_{j=0}^{\infty}\frac{\gamma^{j} r! \left(- q \right)^{\frac{l-r+j}{2}} H_{l-r+j}\left(\frac{ \beta }{2(-q)^{\frac{1}{2}}} \right) }{(r-j)! j!(l-r+j)!}\mathcal{H}_{l}(z_1,z_2|\tau) t^{l-r}, \quad l\geq r, \nonumber
	\end{align}
	where $\xi(\tau):=\frac{\log (1+\tau)}{\tau}$; \quad  $r=0,1,2,\ldots$; \quad $q,\alpha,\beta,\gamma,\delta \in \mathbb{C}$.\\
	
\noindent	\begin{proof}
		 The actions of the 1-parameter groups generated by operators $J^+, J^-, J^3, E, Q$ given by \eqref{hpeq15} for $f \in \mathcal{F}$ are computed in accordance with result \cite[p.18(Theorem 1.10)]{Miller}, so that, we have 
	\begin{align} \label{hpeq26}
		[U(\exp~\beta\mathsf{j}^+&)f](z_1,z_2|\tau;t)=\exp \left(\beta J^+\right)f(z_1,z_2|\tau)=\exp\left(\frac{\log(1+\tau)}{\tau}(\beta^2z_2t^2+\beta z_1t)\right)\nonumber \\
		&\hspace{6cm}\times f(z_1+2z_2\beta t,z_2|\tau;t), \nonumber\\
		[U(\exp~\gamma\mathsf{j}^-&)f](z_1,z_2|\tau;t)=\exp \left(\gamma J^-\right)f(z_1,z_2|\tau)=f\left(z_1+\frac{\gamma}{\alpha t},z_2|\tau;t\right), \nonumber\\
		[U(\exp~\delta\mathsf{j}^3&)f](z_1,z_2|\tau;t)=\exp \left(\delta J^3\right)f(z_1,z_2|\tau)=f(z_1,z_2|\tau;te^\delta),\\
		[U(\exp~\alpha\mathcal{E}&)f](z_1,z_2|\tau;t)=\exp \left(\alpha E\right)f(z_1,z_2|\tau)=\exp(\alpha)f(z_1,z_2|\tau;t), \nonumber\\
		[U(\exp~q\mathcal{Q}&)f](z_1,z_2|\tau;t)=\exp \left(q Q\right)f(z_1,z_2|\tau)=\left(1-4\frac{\log(1+\tau)}{\tau}z_2qt^2\right)^{-\frac{1}{2}} \nonumber \\
		&\hspace{3cm}\times
		\exp\left(\frac{\log(1+\tau)z_1^2}{4\tau z_2}\left(\left(1-4\frac{\log(1+\tau)}{\tau}z_2qt^2\right)^{-1}-1\right)\right) \nonumber \\
		 &\hspace{1cm} \times f\left(z_1\left(1-4\frac{\log(1+\tau)}{\tau}z_2qt^2\right)^{-\frac{1}{2}},z_2|\tau;t\left(1-4\frac{\log(1+\tau)}{\tau}z_2qt^2\right)^{-\frac{1}{2}}\right).\nonumber
	\end{align} 
	
	\noindent Since $g \in K_5$ has parameters $(q,\alpha,b,\gamma,\delta)$, therefore, it has the following expression in terms of the basis element $\mathbf{j}^+,\mathsf{j}^-,\mathsf{j}^3,\mathcal{E},\mathcal{Q}$ of the Lie algebra $\mathcal{K}_5$:
	\begin{equation*}
		g=(\exp \beta\;\mathsf{j}^+)(\exp \gamma\;\mathsf{j}^-)(\exp \delta\;\mathsf{j}^3)(\exp \alpha\mathcal{E})(\exp q\mathcal{Q}).
	\end{equation*}
	\noindent Consequently, we have 
	\begin{equation} \label{hpeq27}
		U(g)f=U(\exp \beta\;\mathsf{j}^+)U(\exp \gamma\;\mathsf{j}^-)U(\exp \delta\;\mathsf{j}^3)U(\exp \alpha\mathcal{E})U(\exp q\mathcal{Q})f.
	\end{equation}
	\noindent A comprehensive computation yields
	\begin{align} \label{hpeq28}
		&[U(g)f](z_1,z_2|\tau;t)=\left(1-4\xi(\tau)z_2qt^2\right)^{-\frac{r+1}{2}}\exp\left(\alpha+r\delta+\xi (\tau)\left(\beta^2z_2t^2+\beta z_1t\right)+ 
		\frac{\xi(\tau)z_1^2}{4z_2} \right. \nonumber \\
		 &\left. \times \left(\left(1-4\xi(\tau)z_2qt^2\right)^{-1}-1\right)\right) f\left(\left(1-4\xi(\tau)z_2qt^2\right)^{-\frac{1}{2}}\left(z_1+2z_2\beta t+\frac{\gamma}{t}\left(1-4\xi(\tau)z_2qt^2\right)\right),z_2|\tau \right), 
	\end{align}
where $\xi(\tau):=\frac{\log (1+\tau)}{\tau}$. From now on, we use $\xi(\tau)$ to denote $\frac{\log (1+\tau)}{\tau}$.\\
	
\noindent	The operator $U$ forms a representations of $K_5$. The entries of matrix elements $A_lr(g)$ of $U(g)$ with respect to the analytic basis $\{f_r(z_1,z_2|\tau;t) = \mathcal{H}_{r}(z_1,z_2|\tau) t^r\}$ are given by
	\begin{equation} \label{hpeq29}
		[U(g)f_r](z_1,z_2|\tau;t)=\sum_{l=0}^{\infty}A_{lr}(g)f_l(z_1,z_2|\tau;t),\quad r=0,1,2,\ldots.
	\end{equation}
	making use of operators \eqref{hpeq28} the above equation taken the form
	\begin{align} \label{hpeq110} 
		&\left(1-4\xi(\tau)z_2qt^2\right)^{-\frac{r+1}{2}}\exp\left(\alpha+r\delta +\xi (\tau)\left(\beta^2z_2t^2+\beta z_1t\right)+\frac{\xi(\tau)}{4 z_2}z_1^2 
		\left(\left(1-4\xi(\tau)z_2qt^2\right)^{-1}-1\right)\right) \nonumber \\ &\hspace{1cm}\times\mathcal{H}_r\left(\left(1-4\xi(\tau)z_2qt^2\right)^{-\frac{1}{2}} \bigg(z_1+2z_2\beta t+ 
		\dfrac{\gamma}{t}\left(1-4\xi(\tau)z_2qt^2\right)\bigg),
		z_2|\tau\right)\nonumber \\
		&\hspace{3cm}=\sum_{l=0}^{\infty}A_{lr}(g) \mathcal{H}_{l}(z_1,z_2|\tau) t^{l-r},\quad
		r=0,1,2,\ldots; \quad q,\alpha,\beta,\gamma,\delta \in \mathbb{C}.
		\end{align}
	Substituting expression \eqref{hpeq17} of matrix elements $A_{lr}(g)$ in equation \eqref{hpeq110}, we get assertion \eqref{hpeq25}.\\
	
\noindent	Again, using alternate expansion \eqref{hpeq85} of matrix elements in equation \eqref{hpeq110}, assertion \eqref{hpeq160} is proved. 
\end{proof}\\

	\noindent	Further, for suitable selection of parameters, the following results are deduced as consequences of Theorem 2.2:\\
	
\noindent	{\bf Corollary 2.1}
		The following implicit summation formula involving the 2-variable degenerate Hermite polynomials $\mathcal{H}_{r}(z_1,z_2|\tau)$ holds true:
	\begin{align} \label{hpeq30}
		&\hspace{-1.5cm}\exp\left(\xi (\tau)\left(\beta^2z_2t^2+\beta z_1t\right)\right) \mathcal{H}_r\left(z_1+2z_2\beta t+\frac{\gamma}{t},z_2|\tau\right) \\
		&=\sum_{l=0}^{\infty} \gamma^{r-l}{L_l}^{(r-l)}(-\beta\gamma)\mathcal{H}_{l}(z_1,z_2|\tau) t^{l-r},\quad  \alpha,\beta,\gamma,\delta \in \mathbb{C} \quad r=0,1,2,\ldots. \nonumber
	\end{align}\\
 
	\noindent	\begin{proof}
		 Taking $q=0$ in equations \eqref{hpeq110} and denoting the corresponding matrix elements by $B_{lr}(g)$, it follows that
		 \begin{align} \label{hpeq111} 
		 	&\hspace{-3cm}\exp\left(\alpha+r\delta +\xi (\tau)\left(\beta^2z_2t^2+\beta z_1t\right)\right)\mathcal{H}_r\left(\left(z_1+2z_2\beta t+ 
		 	\dfrac{\gamma}{t}\right), 
		 	z_2|\tau\right) \\
		 	&=\sum_{l=0}^{\infty}B_{lr}(g) \mathcal{H}_{l}(z_1,z_2|\tau) t^{l-r}, \quad r=0,1,2,\ldots; \alpha,\beta,\gamma,\delta \in \mathbb{C}, \nonumber
		 \end{align}
		which on using expression \eqref{hpeq66} of matrix elements $B_{lr}(g)$, yields assertion \eqref{hpeq30}.  
	\end{proof}\\
	
	\noindent	{\bf Corollary 2.2}
		 The following implicit summation formula involving the 2-variable degenerate Hermite polynomial $\mathcal{H}_{r}(z_1,z_2|\tau)$ holds true:
	\begin{align} \label{hpeq34}
		&\left(1-4\xi(\tau)z_2qt^2\right)^{-\frac{r+1}{2}}\exp\left(\frac{\xi(\tau)z_1^2}{4z_2}\left(\left(1-4\xi(\tau)z_2qt^2\right)^{-1}-1\right)\right) \nonumber \\ 
		 &\times \mathcal{H}_r\left(\left(1-4\xi(\tau)z_2qt^2\right)^{-\frac{1}{2}}\left(z_1+\frac{\gamma}{t}\left(1-4\xi(\tau)z_2qt^2\right)\right),z_2|\tau\right)\nonumber \\
		&=\sum_{l=0}^{\infty} \sum_{j=0}^{\infty} \gamma^{r-l} r! \frac{(q\gamma^2)^j}{j! (2j+r-l)!(l-2j)}\mathcal{H}_{l}(z_1,z_2|\tau) t^{l-r},
	\quad q,\gamma \in \mathbb{C}; \quad r=0,1,2\ldots.
	\end{align}\\

	\noindent \begin{proof}
		 Taking $\alpha=\beta=\delta=0$ in equation \eqref{hpeq110} and denoting the corresponding matrix elements by $C_{lr}(g)$, we have
		\begin{align} \label{hpeq112} 
			&\hspace{-3cm}\left(1-4\xi(\tau)z_2qt^2\right)^{-\frac{r+1}{2}}\exp\left(\frac{\xi(\tau)}{4 z_2}z_1^2 
			\left(\left(1-4\xi(\tau)z_2qt^2\right)^{-1}-1\right)\right) \nonumber \\
			&\hspace{-2cm}\times \mathcal{H}_r\left(\left(1-4\xi(\tau)z_2qt^2\right)^{-\frac{1}{2}} \bigg(z_1+ 
			\dfrac{\gamma}{t}\left(1-4\xi(\tau)z_2qt^2\right)\bigg),
			z_2|\tau\right) \nonumber \\
			&=\sum_{l=0}^{\infty}A_{lr}(g) \mathcal{H}_{l}(z_1,z_2|\tau) t^{l-r}, \quad  r=0,1,2,\ldots; \quad q, \gamma \in \mathbb{C}.
		\end{align}
		 Making use of expression \eqref{hpeq68} of matrix elements $C_{lr}(g)$ in above equation, assertion \eqref{hpeq34} is proved.  
	\end{proof}\\
	
	
	\noindent	{\bf Corollary 2.3}
		 The following implicit summation formula involving the 2-variable degenerate Hermite polynomials $\mathcal{H}_{r}(z_1,z_2|\tau)$ holds true:
	\begin{align} \label{hpeq32}
		&\left(1-4\xi(\tau)z_2qt^2\right)^{-\frac{r+1}{2}}\exp\left(\xi (\tau)\left(\beta^2z_2t^2+\beta z_1t\right)+\frac{\xi (\tau)z_1^2}{4z_2}\left(\left(1-4\xi(\tau)z_2qt^2\right)^{-1}-1\right)\right) \nonumber \\
		 &\times \mathcal{H}_r\left(\left(1-4\xi(\tau)z_2qt^2\right)^{-\frac{1}{2}}\left(z_1+2z_2\beta t\right),z_2|\tau\right) 
		=\sum_{l=0}^{\infty} \frac{(-q)^{\frac{l-r}{2}}}{(l-r)!} \\ &\times \mathcal{H}_{l-r}\left(\frac{\beta}{2(-q)^{\frac{1}{2}}}\right)\mathcal{H}_{l-r}(z_1,z_2,\tau) t^{l-r}, \quad q,\beta \in \mathbb{C}; \quad r=0,1,2\ldots;\quad l\geq r \geq 0. \nonumber
		\end{align}

\noindent	\begin{proof} For $\alpha=\gamma=\delta=0$, equation \eqref{hpeq110} takes the form
	\begin{align} \label{hpeq113} 
		&\left(1-4\xi(\tau)z_2qt^2\right)^{-\frac{r+1}{2}}\exp\left(\xi (\tau)\left(\beta^2z_2t^2+\beta z_1t\right)+\frac{\xi(\tau)}{4\tau z_2}z_1^2 
		\left(\left(1-4\xi(\tau)z_2qt^2\right)^{-1}-1\right)\right) \nonumber \\
		&\hspace{2cm}\times \mathcal{H}_r\left(\left(1-4\xi(\tau)z_2qt^2\right)^{-\frac{1}{2}} \bigg(z_1+2z_2\beta t\bigg),
		z_2|\tau\right)=\sum_{l=0}^{\infty}D_{lr}(g) \\
		& \hspace{4cm}\times \mathcal{H}_{l}(z_1,z_2|\tau) t^{l-r},\quad r=0,1,2,\ldots; \quad q, \beta \in \mathbb{C}, \nonumber
	\end{align}
	which on using expression \eqref{hpeq69} of matrix elements $D_{lr}(g)$, yields summation formula \eqref{hpeq32}.  
\end{proof}\\


\noindent In the next section, implicit summation formulae for certain polynomials related to the 2VDHP $\mathcal{H}_m(z_1,z_2|\tau)$ are derived.
\section{Examples}
We consider the following examples:\\

\noindent	{\bf Example 3.1} For $z_2=-\frac{1}{2}$, the 2VDHP $\mathcal{H}_m(z_1,z_2,\tau)$ reduce to the degenerate Hermite polynomials $\mathcal{H}_m(z|\tau)$, defined by the following generating relation:
\begin{equation} \label{hpeq36}
	\sum_{m=0}^{\infty}\mathcal{H}_{m}(z_1|\tau) \frac{t^m}{m!}=(1+\tau)^{\dfrac{t(z_1-\frac{t}{2})}{\tau}}.
\end{equation}
Therefore, taking $z_2=-\frac{1}{2}$ in equations \eqref{hpeq30}, \eqref{hpeq34}, \eqref{hpeq32}, the following analogues results for $\mathcal{H}_m(z_1|\tau)$ are obtained:
\begin{equation} \label{hpeq37}
	\exp\left(\xi(\tau)\beta z_1t-\frac{1}{2}\xi(\tau)\beta^2 t^2\right) \mathcal{H}_r\left(z_1-\beta t+\frac{\gamma}{t}|\tau\right)=\sum_{l=0}^{\infty} \gamma^{r-l}{L_l}^{r-l}(-\beta\gamma)\mathcal{H}_{l}(z_1|\tau) t^{l-r},
\end{equation}
\begin{align} \label{hpeq39}
	&\left(1+2\xi(\tau)qt^2\right)^{-\frac{r+1}{2}}\exp\left(\frac{\xi(\tau)z_1^2}{2}\left(1-\left(1+2\xi(\tau)qt^2\right)^{-1}\right)\right) \mathcal{H}_r\left(\left(1+2\xi(\tau)qt^2\right)^{-\frac{1}{2}} \right. \nonumber \\ 
	&\hspace{1cm}\left. \times \left(z_1+\frac{\gamma}{t}\left(1+2\xi(\tau)qt^2\right)\right)|\tau\right)=\sum_{l=0}^{\infty} \sum_{j=-\infty}^{\infty} \gamma^{r-l} r! \frac{(q\gamma^2)^j}{j! (2j+r-l)!(l-2j)}\nonumber \\ &\hspace{3cm}\times \mathcal{H}_{l}(z_1|\tau) t^{l-r}, \quad q,\gamma \in \mathbb{C}; \quad r=0,1,2\ldots,
\end{align}

\begin{align} \label{hpeq38}
	&(1+2\xi(\tau)qt^2)^{-\frac{r+1}{2}}\exp\left(\xi(\tau)\beta z_1t-\frac{1}{2}\xi(\tau)\beta^2 t^2-\frac{\beta z_1^2}{2}\left(\left(1+2\xi(\tau)qt^2\right)^{-1}-1\right)\right) \nonumber \\
	&\hspace{1cm}\times \mathcal{H}_r\left(\left(1+2\xi(\tau)qt^2\right)^{-\frac{1}{2}}\left(z_1-\beta t\right)|\tau\right)=\sum_{l=0}^{\infty} \frac{(-q)^{\frac{l-r}{2}}}{(l-r)!}\mathcal{H}_{l-r}\left(\frac{\beta}{2(-q)^{\frac{1}{2}}}\right) \nonumber \\
	&\hspace{3cm}\times	\mathcal{H}_{l-r}(z_1|\tau) t^{l-r}, \quad q,\beta \in \mathbb{C}; \quad r=0,1,2\ldots; \quad l\geq r \geq 0.
\end{align}

\noindent {\bf Example 3.2} 	Taking $q=-\frac{1}{2}$ in equations \eqref{hpeq34} and \eqref{hpeq32}, the following results are obtained:  
\begin{align} \label{hpeq42}
	&\left(1+2\xi(\tau)z_2t^2\right))^{-\frac{r+1}{2}}\exp\left(\frac{\xi(\tau)z_1^2}{4z_2}\left(\left(1+2\xi(\tau)z_2t^2\right)^{-1}-1\right)\right) \nonumber \\
	&\times \mathcal{H}_r\left(\left(1+2\xi(\tau)z_2t^2\right)^{-\frac{1}{2}}\left(z_1+\frac{\gamma}{t}\left(1+2\xi(\tau)z_2t^2\right)\right),z_2|\tau\right)\nonumber \\
	&=\sum_{l=0}^{\infty} \sum_{j=0}^{\infty} \gamma^{r-l} r! \frac{\left(-\frac{\gamma^2}{2}\right)^j}{j! (2j+r-l)!(l-2j)}\mathcal{H}_{l}(z_1,z_2|\tau) t^{l-r},\quad r=0,1,2\ldots; \quad \gamma \in \mathbb{C},
\end{align}
  
\begin{align} \label{hpeq40}
	&(1+2\xi(\tau)z_2t^2)^{-\frac{r+1}{2}}\exp\left(\xi(\tau)\beta^2z_2t^2+\xi(\tau)\beta z_1t+\frac{\xi(\tau)z_1^2}{4z_2}\left(\left(1+2\xi(\tau)z_2t^2\right)^{-1}-1\right)\right) \nonumber \\
	&\hspace{1cm}\times \mathcal{H}_r\left(\left(1+2\xi(\tau)z_2t^2\right)^{-\frac{1}{2}}\left(z_1+2z_2\beta t\right),z_2|\tau\right)=\sum_{l=0}^{\infty} \frac{(2)^{\frac{r-l}{2}}}{(l-r)!}{H}_{l-r}\left(\frac{\beta}{(2)^{\frac{1}{2}}}\right) \nonumber \\
	&\hspace{3cm}\times	\mathcal{H}_{l-r}(z_1,z_2,\tau) t^{l-r}, \quad \beta \in \mathbb{C}; \quad r=0,1,2\ldots;\quad l\geq r \geq 0.
\end{align}
Further, taking $r=0$ in formula \eqref{hpeq40}, we find 
\begin{align} \label{hpeq41}
	&(1+2\xi(\tau)z_2t^2)^{-\frac{1}{2}}\exp\left(\xi(\tau)\beta^2z_2t^2+\xi(\tau)\beta z_1t+\frac{vz_1^2}{4z_2}\left(\left(1+2vz_2t^2\right)^{-1}-1\right)\right) \nonumber \\
	&\hspace{2cm}=\sum_{l=0}^{\infty} \frac{(2)^{\frac{-l}{2}}}{(l)!}{H}_{l}\left(\frac{\beta}{(2)^{\frac{1}{2}}}\right) 
	\mathcal{H}_{l}(z_1,z_2|\tau) t^{l}, \quad \beta \in \mathbb{C};\quad l\geq  0,
\end{align}
which is identical to the Mehler's theorem \cite[p. 140(4.158)]{Miller}.\\

\noindent	{\bf Example 3.3} 
\noindent	Taking $\tau \rightarrow 0$ in equations \eqref{hpeq30}, \eqref{hpeq34}, \eqref{hpeq32}, we find the following summation formulae for $\mathcal{H}_{m}(z_1,z_2)$:
\begin{equation} \label{hpeq43}
	\exp(\beta^2z_2t^2+\beta z_1t) \mathcal{H}_r\left(z_1+2z_2\beta t+\frac{\gamma}{t},z_2\right)=\sum_{l=0}^{\infty} \gamma^{r-l}{L_l}^{r-l}(-\beta\gamma)\mathcal{H}_{l}(z_1,z_2) t^{l-r}, 
	\end{equation}
\begin{align} \label{hpeq45}
	&(1-4z_2qt^2)^{-\frac{r+1}{2}}\exp\left(\frac{z_1^2}{4z_2}\left(\left(1-4z_2qt^2\right)^{-1}-1\right)\right)  \mathcal{H}_m\left(\left(1-4z_2qt^2\right)^{-\frac{1}{2}}\left(z_1+\frac{\gamma}{t}\left(1-4z_2qt^2\right)\right),z_2\right)\nonumber \\
	&=\sum_{l=0}^{\infty} \sum_{j=0}^{\infty} \gamma^{r-l} r! \frac{(q\gamma^2)^j}{j! (2j+r-l)!(l-2j)}\mathcal{H}_{l}(z_1,z_2) t^{l-r},\quad q,\gamma \in \mathbb{C};\;r=0,1,2\ldots,
\end{align}

\begin{align} \label{hpeq44}
	&\hspace{-2cm}(1-4z_2qt^2)^{-\frac{r+1}{2}}\exp\left(\beta^2z_2t^2+\beta z_1t+\frac{z_1^2}{4z_2}\left(\left(1-4z_2qt^2\right)^{-1}-1\right)\right) \nonumber \\ 
	&\hspace{-1cm}\times \mathcal{H}_m\left(\left(1-4z_2qt^2\right)^{-\frac{1}{2}}\left(z_1+2z_2\beta t\right),z_2\right)=\sum_{l=0}^{\infty} \frac{(-q)^{\frac{l-r}{2}}}{(l-r)!}\mathcal{H}_{l-r}\left(\frac{\beta}{2(-q)^{\frac{1}{2}}}\right) \nonumber \\
	&\hspace{1cm}\times	\mathcal{H}_{l-r}(z_1,z_2) t^{l-r}, \quad q,\beta \in \mathbb{C}; \quad r=0,1,2\ldots; \quad l\geq r \geq 0.
\end{align}\\
\noindent In the next section, the Volterra integral equations for the 2VDHP $\mathcal{H}_{m}(z_1,z_2|\tau)$ and other related polynomials are derived. 
\section{Integral Equations}
The connections between Hermite polynomials and other special functions, such as Laguerre polynomials, are crucial for quantum mechanics and mathematical analysis. Hermite polynomials can be introduced in various ways, including generating functions, linear differential equations, power series coefficients, and Volterra integral equations. Integral equations are particularly significant because they ensure uniqueness and offer an analytical approach to physical contexts. Here, we derive the Volterra integral equation for the 2VDHP $\mathcal{H}_{m}(z_1,z_2|\tau)$. Taking $z_1=0$ in equation \eqref{hpeq5}, we have
\begin{equation} \label{hpeq46}
	\mathcal{H}_m(0,z_2|\tau)=\begin{cases}
		\left(\frac{log(1+\tau)}{\tau}\right)^r \frac{(2r)!}{r!}z_2^r &   ,\quad if~ m=2r,\\
		0 & ,\quad if ~m=2r+1,
	\end{cases}
\end{equation}
Similarly, for $z_2=0$, equation \eqref{hpeq5} gives
\begin{equation} \label{hpeq47}
	\mathcal{H}_m(z_1,0|\tau)=
	\left(\frac{log(1+\tau)}{\tau}\right)^m z_1^m.
\end{equation}
Also, from differential relation \eqref{hpeq9}, it follows that
\begin{equation} \label{hpeq48}
	\dfrac{\partial}{\partial z_1}\mathcal{H}_m(z_1,z_2|\tau)\mid_{z_1=0}=\begin{cases}
		\left(\frac{log(1+\tau)}{\tau}\right)^{r+1} \frac{(2r+1)!}{r!}z_2^r & , \quad   if~ m=2r+1,\\
		0 & , \quad if ~m=2r.
	\end{cases}
\end{equation}
Rewriting differential equation \eqref{hpeq14} in the following form:
\begin{equation} \label{hpeq49}
	\left(\dfrac{\partial^2}{\partial z_1^2} + \frac{z_1}{2z_2}\frac{log(1+\tau)}{\tau}\dfrac{\partial}{\partial z_1} -\frac{m}{2z_2}\frac{log(1+\tau)}{\tau}\right)	\mathcal{H}_{m}(z_1,z_2|\tau)=0, \quad m \geq 0.
\end{equation}
 In order to derive the integral equation from differential equation \eqref{hpeq49} subject to initial conditions \eqref{hpeq46} and \eqref{hpeq48}. First, let us focus on the scenario, where $m$ is even, that is for $m=2r$, so that equation \eqref{hpeq49} becomes:
\begin{equation} \label{hpeq50}
	\left(\dfrac{\partial^2}{\partial z_1^2} +\frac{1}{2} \frac{z_1}{z_2}\frac{log(1+\tau)}{\tau}\dfrac{\partial}{\partial z_1} -\frac{2r}{2z_2}\frac{log(1+\tau)}{\tau}\right)	\mathcal{H}_{2r}(z_1,z_2|\tau)=0, \quad r\geq 0.
\end{equation} 
From equation \eqref{hpeq10}, we have
\begin{equation} \label{hpeq51}
	\dfrac{\partial^2}{\partial z_1^2}\mathcal{H}_{2r}(z_1,z_2|\tau)=(2r)(2r-1) \left(\frac{log(1+\tau)}{\tau}\right)^2 \mathcal{H}_{2r-2}(z_1,z_2|\tau).
\end{equation}
Integrating equation \eqref{hpeq51} between the limits $0$ and $z_1$ and using initial condition \eqref{hpeq48}, it follows that
\begin{equation} \label{hpeq52}
	\dfrac{\partial}{\partial z_1}\mathcal{H}_{2r}(z_1,z_2|\tau)=(2r)(2r-1) \left(\frac{log(1+\tau)}{\tau}\right)^2 \int_{0}^{z_1}\mathcal{H}_{2r-2}(\eta,z_2|\tau)d\eta.
\end{equation}
Again, integrating \eqref{hpeq52} and using initial condition \eqref{hpeq46}, it follows that
\begin{align} \label{hpeq53}
	\mathcal{H}_{2r}(z_1,z_2|\tau)=(2r)(2r-1) \left(\frac{log(1+\tau)}{\tau}\right)^2 \int_{0}^{z_1}(z_1-\eta)&\mathcal{H}_{2r-2}(\eta,z_2|\tau)d\eta \nonumber \\ &+\left(\frac{log(1+\tau)}{\tau}\right)^r \frac{(2r)!}{r!}z_2^r.
\end{align}
In view of equations \eqref{hpeq51}$-$\eqref{hpeq53}, differential equation \eqref{hpeq50} takes the form
\begin{align} \label{hpeq54}
	2z_2\mathcal{H}_{2r-2}(z_1,z_2|\tau)- \frac{log(1+\tau)}{\tau} \int_{0}^{z_1}(z_1-\eta)&\{(2r-1)z_1-2r\eta\}\mathcal{H}_{2r-2}(\eta,z_2|\tau)d\eta \nonumber \\
	&-2\left(\frac{log(1+\tau)}{\tau}\right)^{r-1} \frac{(2r-2)!}{(r-1)!}z_2^r=0,
\end{align}
which on replacing $r$ by $r+1$ becomes
\begin{align} \label{hpeq55}
	2z_2\mathcal{H}_{2r}(z_1,z_2|\tau)- \frac{log(1+\tau)}{\tau} \int_{0}^{z_1}(z_1-\eta)&\{(2r+1)z_1-(2r+2)\eta\}\mathcal{H}_{2r}(\eta,z_2|\tau)d\eta \nonumber \\
	&-2\left(\frac{log(1+\tau)}{\tau}\right)^r \frac{(2r)!}{r!}z_2^{r+1}=0.
\end{align}
Moving forward, we examine differential equation \eqref{hpeq49}, for $m=2r+1$, so that we have
\begin{equation} \label{hpeq56}
	\left(\dfrac{\partial^2}{\partial z_1^2} + \frac{z_1}{2z_2}\frac{log(1+\tau)}{\tau}\dfrac{\partial}{\partial z_1} -\frac{2r+1}{2z_2}\frac{log(1+\tau)}{\tau}\right)	\mathcal{H}_{2r+1}(z_1,z_2|\tau)=0, \quad r\geq 0.
\end{equation} 
Using the same arguments as above for differential equation \eqref{hpeq56}, we have
\begin{align} \label{hpeq57}
	2z_2\mathcal{H}_{2r-1}(z_1,z_2|\tau)- \frac{log(1+\tau)}{\tau} \int_{0}^{z_1}(z_1-\eta)&\{2rz_1-(2r+1)\eta\}\mathcal{H}_{2r-1}(\eta,z_2|\tau)d\eta \nonumber \\
	&-2\left(\frac{log(1+\tau)}{\tau}\right)^{r-1} \frac{(2r-1)!}{(r-1)!} z_1z_2^{r}=0.
\end{align}
Again, replacing $r$ by $r+1$ in equation \eqref{hpeq57}, we find 
\begin{align} \label{hpeq58}
	2z_2\mathcal{H}_{2r+1}(z_1,z_2|\tau)- \frac{log(1+\tau)}{\tau} \int_{0}^{z_1}(z_1-\eta)&\{(2r+2)z_1-(2r+3)\eta\}\mathcal{H}_{2r+1}(\eta,z_2|\tau)d\eta \nonumber \\
	&-2\left(\frac{log(1+\tau)}{\tau}\right)^{r} \frac{(2r+1)!}{r!} z_1z_2^{r+1}=0.
\end{align}
Finally, combining equations \eqref{hpeq55} and \eqref{hpeq58}, the following Volterra integral equation for the 2VDHP $\mathcal{H}_m(z_1,z_2|\tau)$ is obtained:
\begin{align} \label{hpeq59}
	2z_2\mathcal{H}_m(z_1,z_2|\tau)- \frac{log(1+\tau)}{\tau} \int_{0}^{z_1}(z_1-\eta)&\{(m+1)z_1-(m+2)\eta\}\mathcal{H}_{m}(\eta,z_2|\tau)d\eta \nonumber \\
	&-2\left(\frac{log(1+\tau)}{\tau}\right)^{n} \frac{(m)!}{n!} z_1^{m-2n} z_2^{n+1}=0,
\end{align}
where $n:=\left[\frac{m}{2}\right]$.\\

\noindent For $z_2=-\frac{1}{2}$, equation \eqref{hpeq5} reduces to generating function \eqref{hpeq36}. Therefore, taking $z_2=-\frac{1}{2}$ and replacing $z_1$ by $z$ in equation \eqref{hpeq59}, the following integral equations for $\mathcal{H}_m(z|\tau)$ is obtained:
\begin{align} \label{hpeq71}
	\mathcal{H}_m(z_1|\tau)+\frac{log(1+\tau)}{\tau} \int_{0}^{z_1}(z_1-\eta)&\{(m+1)z_1-(m+2)\eta\}\mathcal{H}_{m}(\eta|\tau)d\eta \nonumber \\
	&-\left(\frac{-log(1+\tau)}{2\tau}\right)^{n} \frac{(m)!}{n!} z_1^{m-2n}=0,
\end{align}
where $n:=\left[\frac{m}{2}\right]$.\\

\noindent We contemplate the 2-variable degenerate Hermite polynomials $\mathcal{H}_m(z_1,z_2|\tau)$ in this study using the group representation formalism. These polynomials appeared as basis functions for the representation $\uparrow'_{0,1}$ of the $\mathcal{K}_5$. The study in this work validates the feasibility of expanding this methodology to additional practical forms of degenerate Hermite polynomials, which emerged in mathematical physics, applied mathematics, and other scientific fields. 

\end{document}